\documentclass[a4paper,10pt]{article}
\usepackage{geometry,times}
\usepackage{amsfonts}
\usepackage{amssymb}
\usepackage{amsmath}
\usepackage{tikz}
\usepackage{pgfplots}
\usepackage{tikz-3dplot}
\tdplotsetmaincoords{60}{115}
\pgfplotsset{compat=newest}
\definecolor{paleyellow}{rgb}{1, 1, 0.9}
\usepackage{amsthm}

\newcommand{\bb}{\mathbb}

\theoremstyle{plain}
\newtheorem{satz}{prop}
\newtheorem{corollary}[satz]{Corollary}

\newtheorem{definition}[satz]{Definition}

\theoremstyle{thm}
\newtheorem{theorem}[satz]{Theorem}
\theoremstyle{note}

\newtheorem{remark}[satz]{Remark}
\theoremstyle{problem}

\theoremstyle{notadef}
\newtheorem{notation}[satz]{Notation}
\usepackage[]{graphicx}
\usepackage{color}
\definecolor{dark-BLGR}{rgb}{0.28,0.28,0.28}

\definecolor{dark-red}{rgb}{0.90,0.14,0.14}

\definecolor{dark-viol}{rgb}{0.73,0.18,0.69}
\definecolor{Cyan}{rgb}{0.31,0.67,0.82}

\definecolor{light-yellow}{rgb}{1,1,0.8}

\definecolor{Blue}{rgb}{0.0353,0.0275,0.4}

\definecolor{GRAY}{rgb}{0.26,0.26,0.26}

\definecolor{GREYY}{rgb}{0.45,0.45,0.45}

\definecolor{vivid-viol}{rgb}{0.3255,0.0353,0.55}
\newcommand{\Mulberry}{\textcolor{vivid-viol}}
\definecolor{dark-blue}{rgb}{0.05,0.05,0.65}

\definecolor{dark-green}{rgb}{0.03,0.77,0.29}

\definecolor{dark-Green}{rgb}{0.03,0.57,0.09}

\definecolor{strong-viol}{rgb}{0.2353,0.094,0.349}

\definecolor{BLUE}{rgb}{0.41,0.44,0.93}

\definecolor{RED}{rgb}{0.90,0.18,0.32}

\definecolor{Yel}{rgb}{0.89,0.95,0.19}

\definecolor{White}{rgb}{1,1,1}
\definecolor{black}{rgb}{0,0,0}

\definecolor{Orchid}{rgb}{0.6,0.1607,0.2}

\definecolor{Orange}{rgb}{0.99,0.521,0.34902}

\definecolor{magenta}{rgb}{0.99,0.02,0.99}

\definecolor{Thistle}{rgb}{0.1151,0.1249,0.9873}

\definecolor{brown}{rgb}{0.5450,0.2825,0.07453}

 \usepackage{titlesec}
\titleformat{\section}
{\normalfont
  \bfseries
}
{\thesection}{1ex}{}
\titleformat{\subsection}
{\normalfont
  \bfseries
}
    {\thesubsection}{1ex}{}

\newcommand{\NN}{\ensuremath{\mathbb{N}}}

\begin{document}
\hspace*{1cm}\\[.4cm]
{\large \textbf{Explicit formulas for gradients and the divergence in
  $n$-dimensional spherical coordinates}}
\\[.5cm]
Bernd Rummler\\
\hspace*{2mm}{\small \em Otto-von-Guericke-Universit\"at Magdeburg,
  Inst. f\"ur Analysis und Numerik, PF
  4120, 39016 Magdeburg}
\\[2mm]
Gudrun Th\"ater\footnote{Corresponding author\quad
  E-mail:~\textsf{gudrun.thaeter@kit.edu}}\\
\hspace*{2mm}{\small \em Inst. f\"ur Angewandte \& Numerische Mathematik, KIT, 
  76128 Karlsruhe}
\\[5mm]
\hspace*{.1cm}\hfill\parbox{16cm}{
 \\[2mm]
{\small\textbf{Abstract:}
We use the Laplacian in $n$-dimensional spherical coordinates
($n\,\in\,\{2,..,N\}$ with $N\,<\,\infty$) to write the divergence of
a vector field defined on radially symmetric domains ${\Omega}_{(n)}$
in ${\bb R}^{n}$ in the context of vector calculus. 
We apply straightforward equations of  vector calculus with the nabla operator and the transformation matrices from Cartesian to spherical polar coordinates.
One needs the divergence of a vector field e.g. to prove that vector
fields are eigenfunctions 
of the Stokes operator on $n$-dimensional annuli and balls for all
$n\,\in\,{\bb N},\,n>3$, \cite{RuRuTh2026}.
Our divergence formula in partial derivatives in
$n$-dimensional spherical polar coordinates is an important step in a
future 
verification of further Stokes eigenfunctions on those domains. 
}}
\newcommand{\abs}[1]{\left | #1 \right |}
\newcommand{\RR}{\mathbb{R}}
\newcommand{\cC}{{\mathcal{C}}}
\section{Introduction\label{sec_int}}
\begin{figure}[ht]
\hspace*{2cm}
\begin{minipage}[t]{0.5\textwidth}
\begin{tikzpicture}[scale=1.6,tdplot_main_coords] 
\coordinate (P) at ({0},{1.905/sqrt(2)},{1.905/sqrt(2)
});
\shade[ball color = yellow,
    opacity = 0.35
] (0,0,0) circle (1.5cm);
 
\draw[thick,vivid-viol, -stealth] (0,0,0) -- (P);
\draw(0.5,1.5,1.2) node[anchor=east,vivid-viol]{$R$};
%
\end{tikzpicture}
\end{minipage}
\hspace*{-1cm}
\begin{minipage}[t]{0.5\textwidth}
\begin{tikzpicture}[scale=1.6,tdplot_main_coords] 
\coordinate (P) at ({0},{1.905/sqrt(2)},{1.905/sqrt(2)
});
\coordinate (P1) at ({0},{-0.765/sqrt(2)},{0.765/sqrt(2)});
\shade[ball color = yellow,
    opacity = 0.2
] (0,0,0) circle (1.5cm);the the 
\shade[ball color = white,
    opacity = 0.4
] (0,0,0) circle (.75cm);
 
 
 
\draw[thick,vivid-viol, -stealth] (0,0,0) -- (P);
\draw(0.5,1.5,1.2) node[anchor=east,vivid-viol]{$R_o$};
\draw[thick,gray, -stealth] (0,0,0) -- (P1);
\draw((0.25,-0.05,0.35) node[anchor=east]{$R_i$};
%
\end{tikzpicture}
\end{minipage}
\caption{$2$d Sketch of radially symmetric domains ${\Omega}_{(n)}$: ball and annuli with radii
{\Mulberry{$R$\,;}} $R_i$ and $R_o$ (gap-width $R_o-R_i$ ) 
}
\label{fig:geometry}
\end{figure}
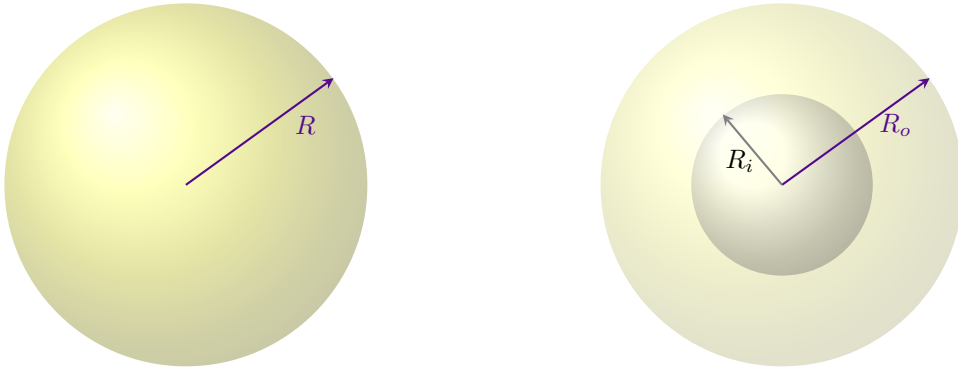

Regarding the flow in radially symmetric $n$-dimensional
domains ${\Omega}_{(n)}$ 
(for two examples see Fig.~\ref{fig:geometry})  one observes,
that the corresponding first Stokes eigenfunction
fixes the related Poincar\'e constant for vector functions on 
domains ${\Omega}_{(n)}$ with vanishing Dirichlet traces.
We have published a series of papers (\cite{RuRuTh2016},
\cite{RuRuTh2025}, \cite{RuRuTh2026}) in this context, 
where we have identified special tools for calculating 
Stokes eigenfunctions. For any finite space dimension see, e.g., \cite{RuRuTh2026}.
There we start with the vector functions 
${\underline{v}}\,=\,{\underline{v}}_{\mathfrak{s}}$ written in
spherical  
polar coordinates to check that the homogeneous Dirichlet boundary
conditions on $\partial{\Omega}_{(n)}$ are fulfilled and that the
divergence of the vector functions ${\underline{v}}$ vanishes. 
The second step then is to verify that every component of
${\underline{v}}\,=\,{\underline{v}}_{\mathfrak{c}}$  
in Cartesian coordinates is an eigenfunction of the Laplacian written in spherical 
polar coordinates. This in order 
to avoid to choose spherical polar coordinates for both, the (vector)
Laplacian $\Delta$ 
and  ${\underline{v}}\,=\,{\underline{v}}_{{\;\!}\mathfrak{s}}$, since
the vector Laplacian when   
applied to ${\underline{v}}_{{\;\!}\mathfrak{s}}\,$ produces
convoluted tensor fields (see \cite{MoonSpencer}), whereas the vector
Laplace operator (written in spherical polar coordinates) 
acts as a scalar on each component of $
{\underline{v}}_{\mathfrak{c}}$ (witten in Cartesian coordinates),
cf. \cite{RumHab}. 

One needs the divergence of the vector functions ${\underline{v}}$ in
partial derivatives  
in $n$-dimensional spherical polar coordinates
to prove that the divergence of the vector functions ${\underline{v}}$
vanishes, since all the components of the vector functions
${\underline{v}}$ are written as 
functions of polar coordinates. Below in Section \ref{Sec2} we are
going to write
the divergence of a function ${\underline{v}}_{{\;\!}\mathfrak{s}}$ using 
vector calculus. This is also a prerequisite to determine 
toroidal fields in ${\bb R}^{3}$  
(cf. also \cite{RumKug1}, \cite{RumTh2024}  and \cite{RuRuTh2025}) -
especially for fields ${\underline{v}}$ in ${\bb R}^{n}$, $n\,>\,3$.

Let us regard our radially symmetric domains ${\Omega}_{(n)}$ in a
non-dimensional setting, where $R$,  $R_i$ and $R_o$ 
are denoting the radius of a ball resp. the inner and outer radius of
the annulus: $0<R<\infty$, $0<R_i<R_o<\infty$.

%
\noindent{\bf General notation A.} Let ${\bb R}^{n}$ be endowed with the
usual Euclidian norm $\| . \|$. 
Elements of ${\bb R}^{n}$ are denoted by underlined small letters.
We write ${\Omega}_{(n)}:=\{{\underline{x}} \in {\bb R}^{n}:\,
\|{\underline{x}}\|\,
<R\}$ for the open balls or ${\Omega}_{(n)}:=\{{\underline{x}} \in {\bb R}^{n}:\,0<R_i <
\|{\underline{x}}\|\,
<R_o\}$  for annuli with radii 
$R_i$ and $R_o$ and use $\omega_{(n)}\,:=\,\{{\underline{x}}\in {\bb
  R}^{n}:\,\|{\underline{x}}\|=1\}$ for the closed surfaces of the open unit balls.
For $r\,\in\,(0,\infty)$ the closed spherical surfaces around the origin with radius $r$ are
$\omega_{(n),r}\,:=\,\{{\underline{x}}\in {\bb
  R}^{n}:\,\|{\underline{x}}\|=r\} $
and the surface areas of $\omega_{(n)}$ are 
$|\omega_{(n)}|\,=\,{2\pi}^{\frac{n}{2}}/
{\Gamma(\frac{n}{2})}$ ($\forall  n\,\in\,{\bb N}$), where ${\Gamma(.)}$ is the ${\Gamma}$-function.\\[1mm]
{\bf General notation B.} Let ${\Omega}_{(n)}$ stand as shorthand for any of the domains
defined above
and the abbreviation $(.)$ for $({\Omega}_{(n)})$, respectively. 
We consider the usual Lebesgue and Sobolev spaces ${\bb L}_{2}(.)$ and
${\bb W}_{2}^{k}(.)$ 
of scalar functions and 
${\underline{\bb L}}_{{\;\!}{2}}(.)=({\bb L}_{2}(.))^{n}$
and ${\underline{\bb W}}_{{\;\!}{2}}^{k}(.)=({\bb
  W}_{2}^{k}(.))^{n}$ of vector functions. 
The norm in ${\bb L}_{2}(.)$ is denoted by $\| . \|_{2}$,
${\bb W}_{2}^{1}\hspace{-.62cm}{~}^{{~}^{{~}^{o}}}\hspace{.2cm}(.)$ is
the closure of $C_{o}^{\infty}(.)$ 
in ${\bb W}_{2}^{1}(.)$.  
All solenoidal vector functions belonging to 
${\underline{C}}_{{\;\!}{o}}^{\infty}(.)$ form $\underline{\cal V}(.)$. The closures
of $\underline{\cal V}(.)$ in ${\underline{\bb L}}_{{\;\!}{2}}(.)$ 
and
${\underline{\bb W}}_{{\;\!}{2}}^{1}(.)$, respectively, are denoted by
${\underline{\bb H}}(.)$ 
and ${\underline{\bb V}}(.)$, respectively.
\begin{notation}\label{Nalba} 
We call the vectorial differential expression
\begin{equation} \label{nab_def}
{\underline{\nabla}}\,:=\,(
\frac{\displaystyle{\partial}}{\displaystyle{\partial x_{1}}},
\frac{\displaystyle{\partial}}{\displaystyle{\partial x_{2}}},\dots,
\frac{\displaystyle{\partial}}{\displaystyle{\partial x_{n}}})^{T}
\end{equation}
the Nabla operator (here ${\underline{\nabla}}$ is written in canonical 
(Cartesian) coordinates).
\end{notation}
\begin{notation}\label{N3} In ${\bb
  R}^{n},\,n\geq 3$, let the unit vectors in the Cartesian coordinate
system 
be given by ${\underline{\mathfrak{e}}}_{j}\,:=\,
(\delta_{j,1},\delta_{j,2}, \dots, \delta_{j,n})^{T}$  ($\forall\,j=1,2,\dots,n$, with
Kronecker's delta $\delta_{j,k}$). 
The polar coordinates are  $r$,  $\vartheta_{1}$, $\dots$ , $\vartheta_{n-2}$ and $\varphi$ with the corresponding unit vectors
${\underline{\mathfrak{e}}}_{{\;\!}r}$, ${\underline{\mathfrak{e}}}_{{\;\!}\vartheta_{1}}$, $\dots$, ${\underline{\mathfrak{e}}}_{{\;\!}\vartheta_{n-2}}$  and
${\underline{\mathfrak{e}}}_{{\;\!}\varphi}$.
If we denote by $\{{\underline{\mathfrak{e}}}_{r},
{\underline{\mathfrak{e}}}_{\vartheta_{1}}, \dots ,
{\underline{\mathfrak{e}}}_{\vartheta_{n-2}},{\underline{\mathfrak{e}}}_{\varphi}\}$
the system of these unit vectors in spherical polar coordinates then ${\underline{u}}$
is representable in both systems as\\[.1cm]
\hspace*{1.3cm}${\underline{u}}\,=\,\sum_{j=1}^{n} u_{j}{\underline{\mathfrak{e}}}_{j}\,=\, \sum_{j=1}^{n} u_{j,{\mathfrak{c}}}{\underline{\mathfrak{e}}}_{j}\,=\, 
 u_{r}{\underline{\mathfrak{e}}}_{r}+
 \sum_{k=1}^{n-2}
 u_{\vartheta_{k}}{\underline{\mathfrak{e}}}_{\vartheta_{k}}+
 u_{\varphi}{\underline{\mathfrak{e}}}_{\varphi}\,=\, 
 u_{r,{\mathfrak{s}}}{\underline{\mathfrak{e}}}_{r}+
\sum_{k=1}^{n-2}
 u_{\vartheta_{k},{\mathfrak{s}}}{\underline{\mathfrak{e}}}_{\vartheta_{k}}
+
 u_{\varphi,{\mathfrak{s}}}{\underline{\mathfrak{e}}}_{\varphi}$.\\[.1cm] 
The transformation from one coordinate system to the other is 
${\underline{u}}_{\mathfrak{c}}\,=\,
{\underline{\underline{T}}}_{{\mathfrak{c}},{\mathfrak{s}}}{\underline{u}}_{\mathfrak{s}}$
or ${\underline{u}}_{\mathfrak{s}}\,=\,
{\underline{\underline{T}}}^{-1}_{{\mathfrak{c}},{\mathfrak{s}}}
{\underline{u}}_{\mathfrak{c}}\,=\,
{\underline{\underline{T}}}_{{\mathfrak{s}},{\mathfrak{c}}}{\underline{u}}_{\mathfrak{c}}$, respectively
(these use columns of coordinates).
The transformation matrices ${\underline{\underline{T}}}_{{\mathfrak{c}},{\mathfrak{s}}}$ and ${\underline{\underline{T}}}_{{\mathfrak{s}},{\mathfrak{c}}}$ are given in the Appendix. 
\end{notation}
\noindent Our paper is organised as follows:
We collect the theoretical background in Section \ref{Sec2} at first. There we outline
the procedures to construct the Laplacian as well as the Stokes operator
as Friedrichs' extension from the Poisson and the Stokes
problem, respectively. We introduce the {\em Leray-Helmholtz projector} 
$\Upsilon : {\underline{\bb
    L}}_{{\;\!}2}(.)\,\longmapsto \,{\underline{\bb H}}(.)$ and sketch the properties of operators with a pure real
point spectrum. 

Subsequently we show in Section \ref{Sec2} how one can  	
take a reading from the Laplacian and the gradient of a scalar
function $\tilde{v}$ in $n$-dimensional spherical polar coordinates to
get the divergence of the vector functions ${\underline{v}}$ in
partial derivatives in $n$-dimensional spherical polar coordinates.
In the appendix we attach the transformation from one coordinate
system to the other for the 
Laplacian.
\section{Theoretical groundwork}\label{Sec2}
\subsection{Laplace and Stokes operators on radially symmetric domains} \label{sec_theo} 
In the following we take both symbols ${\Omega}$ and $(\cdot)$ as
placeholders (as previously in General notation B).
\begin{definition}\label{D4} {\em The Laplace operator is defined in Cartesian coordinates as
\begin{align*}
{\boldsymbol L^{\circledast}}\,{v} := -\Big(
{\displaystyle{\frac{{\partial}^{2} v}{\partial x_{1}^{2}}}}+{\displaystyle{\frac{{\partial}^{2} v}{\partial x_{2}^{2}}}}+\dots\,
{\displaystyle{\frac{{\partial}^{2} v}{\partial x_{n}^{2}}}}\Big)
=- \Delta{\;\!}_{\underline{x} } v
\hspace*{0.7cm}  \forall\,v\,\in\,D({\boldsymbol
  L^{\circledast}})=C_{o}^{\infty}({\Omega})\,. 
\end{align*}
We denote Friedrichs' extension of ${\boldsymbol L^{\circledast}}$
by ${\boldsymbol 
  L}:={\overline{\boldsymbol L^{\circledast}}}$, where ${\boldsymbol
  L}$ is defined on 
$D({\boldsymbol L})\,:=\,{\bb
  W}_{2}^{1}\hspace{-.62cm}{~}^{{~}^{{~}^{o}}}\hspace{.2cm}({\Omega}) 
\cup {\bb W}_{2}^{2}({\Omega})$.}
\end{definition}
\begin{remark} The range of ${\boldsymbol L}$ is
$R({\boldsymbol L})={\bb L}_{2}({\Omega})$. In this sense we may write:
${\boldsymbol L}=-\Delta{\;\!}_{\underline{x} } : D({\boldsymbol
  L})\,\longmapsto \,{\bb L}_{2}(.)$.
\end{remark}
\noindent We need the Leray-Helmholtz projection $\Upsilon$ to define the Stokes operator. 
$\Upsilon$ is the well-defined
projector of ${\underline{\bb L}}_{{\;\!}2}(.)$ onto its subspace 
${\underline{\bb H}}(.)$
of generalised solenoidal fields with vanishing generalised traces in
the normal direction. 
We note, that 
it is also used in the sense of:
$\Upsilon :{\underline{\bb W}}_{{\;\!}2}^{1}(.)\,\longmapsto
\,{\underline{\bb V}}(.)$\,.
\begin{definition}\label{D5} 
{\em The Stokes operator is defined as
$
{\boldsymbol S^{\circledast}}\,{\underline{v}}:=-
\Delta_{\underline{x} } {\underline{v}}\hspace*{0.5cm}
\forall\,{\underline{v}}\in D({\boldsymbol
  S^{\circledast}})=\underline{\cal V}({\Omega})\, 
$, where ${\underline{v}}\,=\,{\underline{v}}_{{\;\!}\mathfrak{c}}$ is written in Cartesian coordinates
and the vector Laplace operator $\Delta_{\underline{x} }$ acts as a
scalar on each component. 
We denote Friedrichs' extension of ${\boldsymbol S^{\circledast}}$ by 
${\boldsymbol S}:={\overline{\boldsymbol S^{\circledast}}}$, where
${\boldsymbol S}$ is defined on its domain 
$D({\boldsymbol S}):={\underline{\bb S}}_{{\;\!}}^{2}(.)=
{\underline{\bb W}}_{{\;\!}2}^{2}(.)\cap{\underline{\bb V}}(.)$ \,.} 
\end{definition}
\begin{remark} The range of ${\boldsymbol S}$ is $R({\boldsymbol
  S})={\underline{\bb H}}(.) $. In this context one may write
${\boldsymbol S} =-\Upsilon \Delta{\;\!}_{\underline{x}
}:{\underline{\bb S}}_{{\;\!}}^{2}(.) 
\,\longmapsto \,{\underline{\bb H}}(.)$. In Definition \ref{D5} one
can also use the Laplace operator in spherical polar
coordinates $\Delta_{r,\vartheta_{1},\dots\,\vartheta_{n-2},\varphi}$
(cf. Remark \ref{Lapldimn}).
We avoid to
choose spherical polar
coordinates for both, $\Delta$
and  ${\underline{v}}\,=\,{\underline{v}}_{{\;\!}\mathfrak{s}}$, since
the vector Laplacian when  
applied to ${\underline{v}}_{{\;\!}\mathfrak{s}}\,$ produces
convoluted tensor fields in this combination 
(see \cite{MoonSpencer}).
\end{remark}
\noindent We sketch the fundamental properties of both operators (i.e. ${\boldsymbol L}$
as well as ${\boldsymbol S}$) using ${\boldsymbol S}$ as an example.
\begin{theorem} \label{thmstok}
The Stokes operator ${\boldsymbol S}$ is positive and self-adjoint.
Its inverse ${\boldsymbol S}^{-1}$ is injective, self-adjoint and compact.
\end{theorem}
\noindent The proof of Theorem \ref{thmstok} is a simple modification of 
Theorems 4.3 and 4.4 in \cite{CoFoi}. The essential tools are
the Rellich theorem and the Lax-Milgram lemma.
The well-known theorem of Hilbert (see, e.g. \cite{CouHil}) and
regularity results like \cite[Prop.~I.2.2]{Temam}
lead to more precise results, namely:
\begin{corollary}
\label{STOeiFU}
The Stokes operator only has a point spectrum.
 All eigenvalues $\lambda_{j}$
of  ${\boldsymbol S}$ are real and of finite multiplicity.
The associated eigenfunctions 
$\{{\underline{w}}_{j}({\underline{x}})\}_{j=1}^{\infty}$
(counted in multiplicity)
are an orthogonal basis of
${\underline{\bb H}}(.)$ and ${\underline{\bb V}}(.)$, i.e.
\begin{align*}
{{(a)}}&\quad {\boldsymbol 
S}{\underline{w}}_{j}:=\lambda_{j}{\underline{w}}_{j}\quad\mbox{for
}\quad{\underline{w}}_{j}\in D({\boldsymbol S})
\quad\forall\,j\in \NN\\
{ (b)}& \quad
0\,<\lambda_{1}\leq\,\lambda_{2}\,\leq\cdots\leq\,\lambda_{j}
\,\leq\cdots\quad\mbox{and}\quad 
\lim_{j\rightarrow\infty}\lambda_{j}=\infty
\\[-2mm]
{(c)}&\quad
\|{\underline{w}}_{j}\|_{{\underline{\bb 
H}}}\,=1\quad\forall \, j\in \NN\,.
\end{align*}
\end{corollary} 
\noindent Concluding we write the eigenvalue problem for the Stokes operator ${\boldsymbol S}$ (cf. Definition \ref{D5}) 
on ${\Omega}_{(n)}$ 
in a classical sense.
\begin{remark}\label{StokesClass}
We look for solutions ${\underline{v}}\,\in\,{\underline{\bb C}}^2(.)$, $\lambda\,\in {\bb R}$ and $p\,\in\,{{\bb C}}^1(.)$  fulfilling the equations:
We write the eigenvalue problem
for the Stokes operator ${\boldsymbol S}$ on ${\Omega}_{(n)}$ in a classical sense, where 
${\underline{v}}\,\in\,{\underline{\bb C}}^2(.)$ and
${\underline{v}}\,=\,{\underline{v}}_{{\;\!}\mathfrak{c}}$ fulfilling is written in Cartesian coordinates
and the vector Laplace operator $\Delta_{\underline{x} }$ and the Nalba operator ${\underline{\nabla}}$
(cf. Notation \ref{Nalba}).
\begin{align}  \label{Sto_Class}
\begin{array}{rcl}
-\,\Delta_{\underline{x}} \,{\underline{v}} \,+\,{\underline{\nabla}} \, p &\,=\,&\lambda {\underline{v}}\\[.2cm]
\mbox{div}{\;\!}\,{\underline{v}}\,=\,{\underline{\nabla}}^T{\underline{v}} &\,=\,& 0 \\[.2cm]
\mbox{with}\,\,{\underline{v}}_{|\partial{{\Omega}_{(n)}}}&\,=\,&{\underline{0}}
\end{array}
\end{align}
Under higher regularity requirements the first two equations ensure that $p$
is a harmonic function, i.e. 
$\Delta_{\underline{x}} \,p\,=\,0\,$.

It is preferable to check the vanishing divergence for vector functions  ${\underline{v}}$ written in spherical polar
coordinates in spherical polar coordinates
due to the simpler access and for one of the first Stokes
eigenfunctions the structure of toroidal fields with  $v_{r}\,=\,0$  
(or $v_{r}{\;\!}{\underline{\mathfrak{e}}}_{r} \,=\,0 \cdot {\;\!}{\underline{\mathfrak{e}}}_{r} $) in spherical polar
coordinates.  We note that toroidal fields ${\underline{v}}$ are solutions of the first equation for the
Stokes eigenfunctions with the
vector Laplacian at $p\,=\,0$ (cf. for $n=3$ \cite{RumTh2024} and \cite{RuRuTh2025}).
\end{remark}
\subsection{The divergence on radially symmetric domains} \label{sec_div_sp}
In the proof  of Theorem 17 in \cite{RuRuTh2026} we needed the divergence 
of a vector function in $n$-dimensional (spherical) polar coordinates.
Unfortunately and to our surprise we did not find this in the
literature in case  $n\,>\,3$. For this we provide an easy and
straightforward way to get 
the divergence in what follows.
We start with  
the definition of the Laplacian with the nabla operator 
in the vector calculus for a scalar function ${\tilde{v}}$:
\begin{align} \label{Trick} 
\Delta\,{\tilde{v}}\,=\,
\mbox{div}{\;\!}\,\mbox{grad}{\;\!}\,{\tilde{v}}\,=\,
{\underline{\nabla}}^{T} \cdot {\underline{\nabla}} \,{\tilde{v}}\,=\,\mbox{div}{\;\!}\,{\underline{v}}_{{\;\!}\mathfrak{s}}\,.
\end{align}
One has to derive ${\underline{v}}_{{\;\!}\mathfrak{s}}$ as the gradient of a scalar function ${\tilde{v}}$ 
by the use of Remark \ref{JacobiMat} (the Jacobian Matrix),  of the reciprocal from the square roots of the $\{g_{j,j}\}_{j=1}^{n}$ there, of the Remark
\ref{R1}  and 
the application of spherical polar coordinates
(cf. Definition \ref{polarcoord}) in a first step.
\begin{remark}[Gradient]\label{Grad_dimn}
For $n\,\geq\,3$  the gradient in spherical coordinates is 
given as
\begin{align*}  
{\underline{v}}_{{\;\!}\mathfrak{s}}\,:=\,
{\underline{\nabla}} \,{\tilde{v}}
\, = \,
\frac{\partial {\tilde{v}}}{\partial r} \cdot {\underline{\mathfrak{e}}}_{{\;\!}r}
+\frac{1}{r}
\frac{\partial {\tilde{v}}}{\partial \vartheta_{1}} \cdot {\underline{\mathfrak{e}}}_{{\;\!}\vartheta_{1}}\,+\,
\frac{1}{r\sin{\vartheta_{1}}}
\frac{\partial {\tilde{v}}}{\partial \vartheta_{2}} \cdot {\underline{\mathfrak{e}}}_{{\;\!}\vartheta_{2}}\,+\,
\frac{1}{r\sin{\vartheta_{1}}\sin{\vartheta_{2}} }
\frac{\partial {\tilde{v}}}{\partial \vartheta_{3}} \cdot {\underline{\mathfrak{e}}}_{{\;\!}\vartheta_{3}}
\,+\,\dots
\quad\\
\,\dots\,
+\frac{1}{r \sin{\vartheta_{1}}\sin{\vartheta_{2}}\cdots  \sin{\vartheta_{n-3} }}
\frac{\partial {\tilde{v}}}{\partial \vartheta_{n-2}} \cdot {\underline{\mathfrak{e}}}_{{\;\!}\vartheta_{n-2}}\,
+\,\frac{1}{r \sin{\vartheta_{1}}\sin{\vartheta_{2}}\cdots  \sin{\vartheta_{n-2} }}
\frac{\partial {\tilde{v}}}{\partial \varphi} \cdot {\underline{\mathfrak{e}}}_{{\;\!}\varphi}\,
\end{align*}
\end{remark}
\begin{remark}[Divergence]\label{Div_dimn} The divergence of 
${\underline{v}}:= {\underline{v}}_{{\;\!}\mathfrak{s}}(r,\vartheta_{1},\dots,\vartheta_{n-2},\varphi)$ in n-dimensional (spherical) polar coordinates is 
\begin{equation} \label{DIV=n2} 
		\mbox{div}{\;\!}\,{\underline{v}}\,=\,
		\frac{1}{r}\left(\frac{\partial{(r\cdot v_{{\;\!}r})}}
		{\partial r}\,+
		\,\frac{\partial{\,v_{{\;\!}\varphi} }}{\partial \varphi}\right)\,\quad\,\mbox{for} \quad \,
n\,=\,2\,\,,
\end{equation}
\begin{equation} \label{DIV=n3} 
		\mbox{div}{\;\!}\,{\underline{v}}\,=\,
		\frac{1}{r^2}\frac{\partial{(r^2\cdot v_{{\;\!}r})}}
		{\partial r}\,+\,
		\frac{1}{r \sin\vartheta_{1}}\left(
		\frac{\partial(\,{\sin\vartheta_{1}}\,v_{{\;\!}\vartheta_{1}})}{\partial \vartheta_1}\,+\,
		\frac{\partial{\,v_{{\;\!}\varphi} }}{\partial \varphi}\right)\,
		\,\quad\,\mbox{for} \quad \,
n\,=\,3\,\,,
\end{equation}
\begin{equation}\label{DIV=n4}
\mbox{div}{\;\!}\,{\underline{v}}\,=\,
		\frac{1}{r^3}\frac{\partial{(r^3\cdot v_{{\;\!}r})}}
		{\partial r}\,+\,
		\frac{1}{r \sin^2\vartheta_{1}}
		\frac{\partial(\,{\sin^2\vartheta_{1}}\,v_{{\;\!}\vartheta_{1}})}{\partial \vartheta_1}
		\,+\,\frac{1}{r \sin\vartheta_{1}\sin{\vartheta_{2}}}\left(
		\frac{\partial(\,{\sin\vartheta_{2}}\,v_{{\;\!}\vartheta_{2}})}{\partial \vartheta_2}\,+\,
\frac{\partial{\,v_{{\;\!}\varphi} }}{\partial \varphi}\right)\,		
\quad\,\mbox{for} \quad \,
n\,=\,4\,\,
\end{equation}
and $\quad \dots$ \,for \,\,$n\,=\,n$ :
\begin{align}\label{DIV=nn}
\hspace*{-1.3cm}
\mbox{div}{\;\!}\,{\underline{v}}\,=\,\frac{1}{r^{n-1}}\frac{\partial{(r^{n-1}\cdot v_{{\;\!}r})}}
		{\partial r}\,+\,\frac{1}{r \sin^{n-2}\vartheta_{1}}
		\frac{\partial(\,{\sin^{n-2}\vartheta_{1}}\,v_{{\;\!}\vartheta_{1}})}{\partial \vartheta_1}  
		\,+\,
		\frac{1}{r \sin\vartheta_{1} {\sin^{n-3}\vartheta_{2}}}
		\frac{\partial(\,{\sin^{n-3}\vartheta_{2}}\,v_{{\;\!}\vartheta_{2}})}{\partial \vartheta_{2}}\,+\,\dots\,+
		\quad\quad\\
+\,\frac{1}{r \sin{\vartheta_{1}}\sin{\vartheta_{2}}\cdots \sin{\vartheta_{n-4}}\sin^2{\vartheta_{n-3}} }
\frac{\partial(\,{\sin^2\vartheta_{n-3}}\,v_{{\;\!}\vartheta_{n-3}})}{\partial \vartheta_{n-3}}\,+\,
\frac{1}{r \sin{\vartheta_{1}}\sin{\vartheta_{2}}\cdots \sin{\vartheta_{n-2}}}\left(
		\frac{\partial(\,{\sin\vartheta_{n-2}}\,v_{{\;\!}\vartheta_{n-2}})}{\partial \vartheta_{n-2}}\,+\,\frac{\partial{\,v_{{\;\!}\varphi} }}{\partial \varphi}\right)\,,	\nonumber 
\,
\end{align}
where the above statements follow from straightforward standard calculations.
\end{remark}

\section*{Appendix}\label{App_Bes}
Let the unit vectors in the Cartesian coordinate system in ${\bb
  R}^{n},\,n\geq 3$
be given by ${\underline{\mathfrak{e}}}_{j}\,:=\,
(\delta_{j,1},\delta_{j,2}, \dots, \delta_{j,n})^{T}$ for all $j=1,2,\dots,n$, with
Kronecker's delta $\delta_{j,k}$. 
The polar coordinates are stated as $r$,  $\vartheta_{1}$, $\dots$ , $\vartheta_{n-2}$ and $\varphi$ 
(cf. the following Definition \ref{polarcoord}) with the corresponding unit vectors
${\underline{\mathfrak{e}}}_{{\;\!}r}$, ${\underline{\mathfrak{e}}}_{{\;\!}\vartheta_{1}}$, $\dots$, ${\underline{\mathfrak{e}}}_{{\;\!}\vartheta_{n-2}}$  and
${\underline{\mathfrak{e}}}_{{\;\!}\varphi}$.
\begin{definition} \label{polarcoord}
The representation of any point ${\underline{x}} \in {\bb R}^{n}$ in the system
of polar coordinates is given via
\begin{align}  \label{Pol1Sh}
\begin{array}{rcl}
x_1 &\,=\,&r\cdot \sin{\vartheta_{1}}\cdots\sin{\vartheta_{n-2}}\cos{\varphi}\\
x_2 &\,=\,&r\cdot \sin{\vartheta_{1}}\cdots\sin{\vartheta_{n-2}}\sin{\varphi}\\
x_3 &\,=\,&r\cdot \sin{\vartheta_{1}}\cdots\cos{\vartheta_{n-2}}\\
{~} & \vdots  & {~} \\
x_{n-1} &\,=\,&r\cdot \sin{\vartheta_{1}}\cos{\vartheta_{2}}\\
x_{n} &\,=\,&r\cdot \cos{\vartheta_{1}}
\end{array}\\
\mbox{with}\,\,
r\,:=\,\|{\underline{x}}\| \in [0,\infty) ,\,\vartheta_{1},\dots,\vartheta_{n-2}\,\in\,[0,\pi]
\,\,\mbox{and}\,\,\varphi\,\in\,[0,2\pi] \,\,.\nonumber
\end{align}
\end{definition}
\begin{notation}[Surface harmonics of degrees $\ell = 0$ and $\ell = 1$]\label{SurfHarm0+1}
The function $f(\vartheta_{n-2},\dots\,\vartheta_{1},\varphi)\,=\,const.\,\not=\,0$  is a non-vanishing harmonic 
polynomial of degree $\ell = 0$ in $r$. For any $ n\,\in {\bb
  N}:n\,>\,1$ we state the sherical  
surface harmonic function
\begin{align}  \label{SphSurf0}
S^{\{0\}}\,:=\,1\,.
\end{align}
In the representation \eqref{Pol1Sh} the functions $\{x_{k}\}_{k=1}^n$
are harmonic polynomials
of degree $\ell = 1$ in $r$. We call the functions
 \begin{align}  \label{SphSurf1}
\begin{array}{rcl}
S^{\{1\}}_1 &\,=\,& \sin{\vartheta_{1}}\cdots\sin{\vartheta_{n-2}}\cos{\varphi}\\
S^{\{1\}}_2 &\,=\,& \sin{\vartheta_{1}}\cdots\sin{\vartheta_{n-2}}\sin{\varphi}\\
S^{\{1\}}_3 &\,=\,& \sin{\vartheta_{1}}\cdots\cos{\vartheta_{n-2}}\\
{~} & \vdots  & {~} \\
S^{\{1\}}_{n-1} &\,=\,& \sin{\vartheta_{1}}\cos{\vartheta_{2}}\\
S^{\{1\}}_{n} &\,=\,& \cos{\vartheta_{1}}\,,
\end{array}\\
\mbox{where}\,\,
\,\vartheta_{1},\dots,\vartheta_{n-2}\,\in\,[0,\pi]
\,\,\mbox{and}\,\,\varphi\,\in\,[0,2\pi] \,\,,\nonumber
\end{align}
sherical surface harmonics of degree $\ell = 1$. We write $S\,\in\,{\mbox{span}}\{S^{\{1\}}_k\}_{k=1}^n$ 
for a spherical surface harmonic function of degree $\ell = 1$
(e.g. the Definition in Subsection  6.3.1 in \cite{Triebel}) as well.
\end{notation}
\begin{remark}\label{JacobiMat}
The first step in the calculation of
the transformation between the Cartesian and spherical polar coordinates
(cf. Definition \ref{polarcoord})
is to calculate the Jacobian matrix:
\begin{align*}
{\underline{\underline{J}}} =
\begin{bmatrix} 
\frac{\partial x_1}{\partial r}&\frac{\partial x_1}{\partial \vartheta_{1}}&
\dots &\frac{\partial x_1}{\partial \vartheta_{n-2}}
&\frac{\partial x_1}{\partial{\varphi}}\\
\vdots&\vdots & {~}& \vdots&\vdots\\
\frac{\partial x_n}{\partial r}&\frac{\partial x_n}{\partial \vartheta_{1}}&
\dots  &\frac{\partial x_n}{\partial \vartheta_{n-2}}&\frac{\partial x_n}{\partial{\varphi}}
\end{bmatrix}
  \end{align*}
 and in a second step the corresponding metric tensor 
 ${\underline{\underline{g}}}$:
\begin{align*}  
{\underline{\underline{g}}}\,:=\,{\underline{\underline{J}}} \cdot {\underline{\underline{J}}}^T\,=\,
\text{diag}\{g_{j,j}\}_{j=1}^{n}\quad\,\text{with}\,\, g_{1,1}=1,\,g_{2,2}=r^2,\,\dots\,,
g_{n,n}=r^2{\sin}^2{\vartheta_{1}}  \cdots {\sin}^2{\vartheta_{n-2}}
\,.
\end{align*}
The reciprocal of the square roots of the $\{g_{j,j}\}_{j=1}^{n}$ used
as multipliers column by column 
applied on ${\underline{\underline{J}}}$ provide  the transformation matrices.
\end{remark}
\begin{remark}
\label{R1} The transformation between the Cartesian coordinates and the spherical polar coordinates
(cf. Notation \ref{N3}) as the transformation of one coordinate system to the other 
can be written as ${\underline{u}}_{\mathfrak{c}}\,=\,
{\underline{\underline{T}}}_{{\mathfrak{c}},{\mathfrak{s}}}{\underline{u}}_{\mathfrak{s}}$
\,or\,\,\,${\underline{u}}_{\mathfrak{s}}\,=\,
{\underline{\underline{T}}}^{-1}_{{\mathfrak{c}},{\mathfrak{s}}}
{\underline{u}}_{\mathfrak{c}}\,=\,
{\underline{\underline{T}}}_{{\mathfrak{s}},{\mathfrak{c}}}{\underline{u}}_{\mathfrak{c}}$,
respectively,
where we  use the concept of columns of coordinates
and the transformation matrices ${\underline{\underline{T}}}_{{\mathfrak{s}},{\mathfrak{c}}}\,:=\,
{\underline{\underline{T}}}_{{\mathfrak{c}},{\mathfrak{s}}}^{-1}\,=\,
{\underline{\underline{T}}}_{{\mathfrak{c}},{\mathfrak{s}}}^{T}\,$ 
\begin{equation} \label{R1_T=2} 
		{\underline{\underline{T}}}_{{\mathfrak{c}},{\mathfrak{s}}}\,
		:=\,\left[
		\begin{array}{lc}
\cos{\varphi}& -\sin{\varphi} \\
\sin{\varphi}&  {~}\cos{\varphi}
		\end{array}
		\right] \quad\,\mbox{for} \quad \,
n\,=\,2\,\,,
\end{equation}
\begin{equation}\label{R1_T=3} 
		{\underline{\underline{T}}}_{{\mathfrak{c}},{\mathfrak{s}}}\,
		:=\,\left[
		\begin{array}{llc}
\sin{\vartheta_{1}}\cos{\varphi}& \cos{\vartheta_{1}}\cos{\varphi}& -\sin{\varphi} \\
\sin{\vartheta_{1}}\sin{\varphi}& \cos{\vartheta_{1}}\sin{\varphi}& {~}\cos{\varphi}\\	
\cos{\vartheta_{1}} & -\sin{\vartheta_{1}}& {~} 0
		\end{array}
		\right] \quad\,\mbox{for} \quad \,
n\,=\,3\,\,,
\end{equation}
\begin{equation}\label{R1_T=4} 
		{\underline{\underline{T}}}_{{\mathfrak{c}},{\mathfrak{s}}}\,
		:=\,\left[
		\begin{array}{llcc}
\sin{\vartheta_{1}}\sin{\vartheta_{2}}\cos{\varphi} & \cos{\vartheta_{1}}\sin{\vartheta_{2}}\cos{\varphi}& \cos{\vartheta_{2}}\cos{\varphi}&
-\sin{\varphi} \\
\sin{\vartheta_{1}}\sin{\vartheta_{2}}\sin{\varphi} & \cos{\vartheta_{1}}\sin{\vartheta_{2}}\sin{\varphi}& \cos{\vartheta_{2}}\sin{\varphi} & {~}\cos{\varphi}\\
\sin{\vartheta_{1}}\cos{\vartheta_{2}} & \cos{\vartheta_{1}}\cos{\vartheta_{2}} & -\sin{\vartheta_{2}}& {~} 0\\
\cos{\vartheta_{1}} & -\sin{\vartheta_{1}}& {~} 0 & {~} 0
		\end{array}
		\right] \quad\,\mbox{for} \quad \,
n\,=\,4\,\,,
\end{equation}
and for general $n=n$
\begin{equation}\label{R1_T=n} 
		{\underline{\underline{T}}}_{{\mathfrak{c}},{\mathfrak{s}}}\,
		:=\,\left[
\begin{array}{llccc}
\sin{\vartheta_{1}}\sin{\vartheta_{2}}\cdots \sin{\vartheta_{n-2}}\cos{\varphi} & 
\cos{\vartheta_{1}}\sin{\vartheta_{2}}\cdots \sin{\vartheta_{n-2}}\cos{\varphi}& \dots &
\cos{\vartheta_{n-2}}\cos{\varphi}&
-\sin{\varphi} \\
\sin{\vartheta_{1}}\sin{\vartheta_{2}}\cdots \sin{\vartheta_{n-2}}\sin{\varphi} & 
\cos{\vartheta_{1}}\sin{\vartheta_{2}}\cdots \sin{\vartheta_{n-2}}\sin{\varphi} & \dots &\cos{\vartheta_{n-2}}\sin{\varphi}&
 {~} \cos{\varphi}\\
\sin{\vartheta_{1}}\sin{\vartheta_{2}}\cdots \cos{\vartheta_{n-2}}& 
\cos{\vartheta_{1}}\sin{\vartheta_{2}}\cdots \cos{\vartheta_{n-2}} & \dots & - \sin{\vartheta_{n-2}}&
 {~} 0\\
\vdots & \vdots & \ddots &\vdots & \vdots\\
\sin{\vartheta_{1}}\cos{\vartheta_{2}} & \cos{\vartheta_{1}}\cos{\vartheta_{2}} &\dots &
 {~} 0 & {~} 0\\
\cos{\vartheta_{1}} & -\sin{\vartheta_{1}}& \dots &    {~} 0 & {~} 0
		\end{array}
		\right] 
                \,.
\end{equation}
\end{remark}
\begin{remark}\label{Lapldimn}
The Laplacian in spherical coordinates
$\Delta_{sph}(.):=\Delta_{r,\vartheta_{1},\dots\,\vartheta_{n-2},\varphi}(.)$ is
($n\geq3$; see, e.g.,  6.3.4 in \cite{Triebel})
\begin{align*} \hspace*{-.3cm}
 \Delta_{sph}(.)  = \frac{1}{r^{n-1}}
\frac{\partial{ }}{\partial r}({r^{n-1}}\frac{\partial { (.)}}{\partial r})
+\frac{1}{r^{2}}\left(
\frac{1}{\sin^{n-2}{\vartheta_{1}}}\frac{\partial{ }}{\partial \vartheta_{1}}
(\sin^{n-2}\vartheta_{1}\frac{\partial{ }(.)}{\partial \vartheta_{1}})+
\frac{1}{\sin^{2}\vartheta_{1}\sin^{n-3}{\vartheta_{2}}}\frac{\partial{ }}{\partial \vartheta_{2}}
(\sin^{n-3}\vartheta_{2}\frac{\partial{ }(.)}{\partial \vartheta_{2}})\,+\right.
\\
{~} 
\left.\,\dots\,+\,
\frac{1}{\sin^{2}\vartheta_{1}\sin^{2}{\vartheta_{2}}\cdots
\sin^{2}{\vartheta_{n-3}}
\sin^{~}{\vartheta_{n-2}}
}\frac{\partial{ }}{\partial \vartheta_{n-2}}
(\sin \vartheta_{n-2}\frac{\partial{ }(.)}{\partial \vartheta_{n-2}})\,+
\,
\frac{1}{\sin^{2}\vartheta_{1}\sin^{2}{\vartheta_{2}}\cdots
\sin^{2}{\vartheta_{n-2}}}
  \frac{\partial^{2}{ (.)}}{\partial \varphi^{2}}\right)\,\,,
\end{align*}
resp. $\displaystyle
\Delta_{r,\vartheta_{1},\dots\,\vartheta_{n-2},\varphi}(.)  = 
{~} 
\,\frac{1}{r^{n-1}}
\frac{\partial{ }}{\partial r}({r^{n-1}}\frac{\partial { (.)}}{\partial r})
\,-\,\frac{\displaystyle{1}}{\displaystyle{r^{2}}}{ \mbox{B}}{\;\!}(.)\,,\hspace*{.2cm}{~}
$
where B$(.) $ denotes the Beltrami differential operator.
\end{remark}
\noindent We define the Laplace-Beltrami operator by means of
Beltrami's differential operator
in the following
\begin{definition}\label{DBeltr} { For all $Y\,\in\,D({\boldsymbol
  B^{\circledast}})=C^{\infty}({\omega_{(n)}}) \subset {\bb
  L}_{2}({\omega_{(n)}})\,$  the Laplace-Beltrami operator is defined
as}
\begin{align*}
{\boldsymbol B^{\circledast}}\,{Y} := \,{ \mbox{B}}{\;\!}(Y)\,\,=\,-
\left(
\frac{1}{\sin^{n-2}{\vartheta_{1}}}\frac{\partial{ }}{\partial \vartheta_{1}}
(\sin^{n-2}\vartheta_{1}\frac{\partial{ } Y}{\partial \vartheta_{1}})+
\frac{1}{\sin^{2}\vartheta_{1}\sin^{n-3}{\vartheta_{2}}}\frac{\partial{ }}{\partial \vartheta_{2}}
(\sin^{n-3}\vartheta_{2}\frac{\partial{ } Y }{\partial \vartheta_{2}})\,+\right.\hspace*{3.1cm}  
\\
{~} 
\left.\,\dots\,+\,
\frac{1}{\sin^{2}\vartheta_{1}\sin^{2}{\vartheta_{2}}\cdots
\sin^{2}{\vartheta_{n-3}}
\sin^{~}{\vartheta_{n-2}}
}\frac{\partial{ }}{\partial \vartheta_{n-2}}
(\sin \vartheta_{n-2}\frac{\partial{ } Y}{\partial \vartheta_{n-2}})\,+
\,
\frac{1}{\sin^{2}\vartheta_{1}\sin^{2}{\vartheta_{2}}\cdots
\sin^{2}{\vartheta_{n-2}}}
\frac{\partial^{2}{\, Y }}{\partial \varphi^{2}}\right)\,\,.
\end{align*}
We denote the Friedrichs' extension of ${\boldsymbol B^{\circledast}}$
by ${\boldsymbol 
  B}:={\overline{\boldsymbol B^{\circledast}}}$, where ${\boldsymbol
  B}$ is applied on 
$D({\boldsymbol B})\,:=\,{\bb W}_{2}^{2}({\omega_{(n)}})\subset {\bb L}_{2}({\omega_{(n)}})$.
\end{definition}
\begin{remark} The detailed construction of the Laplace-Beltrami operator ${\boldsymbol
  B}$ is given in \cite[Subsection 6.3.5]{Triebel} at great length. Especially, the step from an $n$-dimensional shell to the boundary ${\omega_{(n)}}$ is illustrated there.
\end{remark}
\noindent We cite explicitly the following result which is important for the eigenfunction of the Laplacian 
as well as for the first Stokes eigenfunctions:
\begin{theorem} \label{thmsbeltr}
The Laplace-Beltrami operator ${\boldsymbol B}$ is nonnegative and self-adjoint.
${\boldsymbol B}$ is an operator with pure point spectrum. Its eigenvalues are $\ell (\ell+n-2)$,  
$\ell\,=\,0,1,2,\dots $.
The surface harmonics S(.) of the degree $\ell$ form a set of all eigenfunctions of ${\boldsymbol B}$
to the eigenvalue $\ell (\ell+n-2)$.
\end{theorem}
\end{document}